\documentclass[11pt,leqno,oneside,letterpaper, openany]{amsart}
\usepackage[letterpaper,left=1.75truein, top=1truein]{geometry}
\usepackage{amsmath,amsfonts,amssymb,amscd,amsthm,amsbsy,epsf,graphics} 
\usepackage{stmaryrd}
\usepackage{setspace}
\usepackage[title]{appendix}
\usepackage{comment} 
\usepackage[normalem]{ulem}

\usepackage[dvipsnames]{xcolor}
\usepackage{color}

\definecolor{org}{RGB}{225,119,0}

\usepackage[colorlinks,linkcolor=blue,citecolor=blue]{hyperref} 

\usepackage{bm}
\usepackage{caption,subcaption,graphicx,epsfig}
\usepackage{tikz-cd,calc}
\usetikzlibrary{knots}
\usetikzlibrary{decorations.markings}

\usepackage{lipsum}

\makeatletter
\def\l@subsection{\@tocline{2}{0pt}{4pc}{5pc}{}}
\makeatother
 
\let\oldtocsection=\tocsection
 
\let\oldtocsubsection=\tocsubsection
 
\let\oldtocsubsubsection=\tocsubsubsection
 
\renewcommand{\tocsection}[2]{\hspace{0em}\oldtocsection{#1}{#2}}
\renewcommand{\tocsubsection}[2]{\hspace{0em}\oldtocsubsection{#1}{#2}}
\renewcommand{\tocsubsubsection}[2]{\hspace{2em}\oldtocsubsubsection{#1}{#2}}

\DeclareGraphicsExtensions{.pdf}

\newtheorem{theorem}{Theorem}[section] 
\newtheorem*{theorem*}{Theorem}
\newtheorem{corollary}[theorem]{Corollary}
\newtheorem*{corollary*}{Corollary}
\newtheorem{lemma}[theorem]{Lemma}
\newtheorem{proposition}[theorem]{Proposition}
\newtheorem*{proposition*}{Proposition}

\newtheorem{question}[theorem]{Question}
\newtheorem*{question*}{Question}

\theoremstyle{definition}

\newtheorem*{problem*}{Problem}

\newtheorem{remark}[theorem]{Remark} 
\newtheorem*{remark*}{Remark}

\newtheorem*{acknowledgement*}{Acknowledgements}

\newtheoremstyle{cases}
  {12pt plus 6 pt}
  {2pt}
  {\bfseries}   
  {}
  {\bfseries}
  {.}
  {.5em}
  {}
\theoremstyle{cases}

\numberwithin{subcase}{case} 
\numberwithin{subsubcase}{subcase}
\numberwithin{equation}{subsection} 

\graphicspath{{figures/}} 

\begin{document}

\title[Recognition problem for 4-manifolds; unknotting problem for 2-spheres]{The recognition problem for 4-manifolds and the unknotting problem for 2-spheres in 4-manifolds}

\author[Cameron McA. Gordon]{Cameron McA. Gordon} 
\address{University of Texas at Austin, Department of Mathematics - PMA 8.100, 2515 Speedway, Stop C1200, Austin, TX 78712, USA}                                                                                 
\email{gordon@math.utexas.edu}



\maketitle

\begin{abstract}
We show that there is no algorithm to decide whether or not a given 4-manifold is PL homeomorphic to the connected sum of 7 copies of $S^2 \times S^2$. We also show that the TOP and PL unknotting problems for all surfaces in the connected sum of 7 copies of $S^2 \times S^2$ are unsolvable.
\end{abstract} 

\section{Introduction}

In this paper we address two algorithmic questions in 4-dimensional topology: the homeomorphism problem for 4-manifolds and the unknotting problem for 2-spheres in 4-manifolds.

To ensure that the objects under consideration can clearly be described by a finite amount of data we assume that our manifolds and knots are PL. Thus the homeomorphism problem asks whether or not two given closed, connected PL manifolds are homeomorphic. There are two versions: we can ask if they are PL homeomorphic or merely homeomorphic (TOP homeomorphic). (We will use CAT to denote either TOP or PL.) In [M] Markov showed that both these problems are algorithmically unsolvable. The proof is based on the existence, proved by Novikov [N] and Boone [B], of a finitely presented group with unsolvable word problem. Given such a group, Adjan [A] and Rabin [R] showed how to construct, for some $k$, a {\it $k$-relator Adjan-Rabin set}, i.e. a recursively enumerable set $\mathcal P$ of finite $k$-relator group presentations such that there is no algorithm to decide whether or not the group presented by a given $P \in \mathcal P$ is trivial. Markov showed that if there exists such a set then the CAT {\it recognition problem} for $\#_k (S^2 \times S^2)$ is unsolvable: there is no algorithm to decide whether or not a given PL 4-manifold is CAT homeomorphic to $\#_k (S^2 \times S^2)$. 

Borisov constructed a 12-relator presentation of a group with unsolvable word problem [Bo], and the author showed that from a finite $n$-relator presentation of a group with unsolvable word problem one could construct an $(n+2)$-relator Adjan-Rabin set [G1]. It follows that the CAT recognition problem for $\#_{14} (S^2 \times S^2)$ is unsolvable; see [St1],[St2],[CL]. In [G2] the author made minor improvements to both the algebraic and topological parts of the argument and showed that the CAT recognition problem for $\#_{12} (S^2 \times S^2)$ is unsolvable.

A new approach was taken by Tancer in [T]. Defining the {\it deficiency} of a group presentation with $m$ generators and $n$ relators to be $n - m$, he showed, using Freedman's topological classification of simply-connected 4-manifolds [F], that if there exists an {\it Adjan-Rabin set of deficiency $d$} (i.e. each presentation has deficiency $d$), then the TOP recognition problem for $\#_d (S^2 \times S^2)$ is unsolvable. Moreover, Tancer constructed an Adjan-Rabin set of deficiency 9, thereby showing that the TOP recognition problem for  $\#_9 (S^2 \times S^2)$ is unsolvable.

In an earlier version of this work we proved the unsolvability of the CAT recognition problem for $\#_8 (S^2 \times S^2)$. To do this we modified Borisov's argument by proving a more general version of the embedding theorem in [Bo, Section 3], which enabled us to construct an Adjan-Rabin set of deficiency 8.

Recently Kegel, Li and Ren made a major improvement to Borisov's construction to get a 3-generator 9-relator presentation of a group with unsolvable word problem [KLR]. Using this they produced an Adjan-Rabin set of deficiency 7, proving the unsolvability of the TOP recognition problem for $\#_7 (S^2 \times S^2)$, and a 9-relator Adjan-Rabin set, proving the unsolvability of the PL recognition problem for $\#_9 (S^2 \times S^2)$. We use their group in this paper. Similarly, for Theorem \ref{theorem: univ} we use the 13-relator universal group given by Ren [Re], which improves on the 21-relator example of Valiev [V].  

We show that Tancer's approach can be made to work in the PL category, giving the following theorem.

\begin{theorem}
  \label{theorem: rec7}
  The PL recognition problem for $\#_7 (S^2 \times S^2)$ is unsolvable.
\end{theorem}

The ultimate goal would be to reduce the number of $S^2 \times S^2$ summands to zero.

\begin{question}
  \label{question: s4rec}
Is the CAT recognition problem for $S^4$ unsolvable?
\end{question}

The CAT recognition problem for $S^n$ is unsolvable for $n \ge 5$ [VKF, Appendix by S.P. Novikov], and solvable for $n = 3$ [Ru],[Th] and $n = 2$. We note that in these dimensions a PL $n$-manifold is TOP homeomorphic to $S^n$ if and only if it PL homeomorphic to $S^n$ [Mo1],[S]. (The last is for the smooth case; the proof can be adapted to work in the PL category.) 

We now turn to codimension 2 unknotting problems. The most basic of these is the unknotting problem for {\it $n$-knots}: is there an algorithm to decide whether or not a given locally flat PL embedded $n$-sphere in $S^{n+2}, n \ge 1$, is unknotted? This is solvable for $n = 1$ [H], unsolvable for $n \ge 3$ [NW],[GGS], and unknown for $n = 2$. Note that for $n$-knots with $n = 1$ or $n \ge 3$, PL and TOP unknotted are equivalent [Mo2],[L1]. Again, this is unknown in the case $n = 2$.

Thus, in this setting the analog of Question \ref{question: s4rec} is the following.

\begin{question}
\label{question: 2knot}
Is the CAT unknotting problem for 2-spheres in $S^4$ solvable?
\end{question}

In the absence of an answer to this question, there are two natural generalizations that one might consider: higher genus surfaces in $S^4$, and 2-spheres in other 4-manifolds. The first is considered in [G3], which for example proves the unsolvability of the TOP and PL unknotting problem for orientable surfaces in $S^4$ of genus $\ge 18$ and $19$, respectively. Here we consider the second generalization, and prove 

\begin{theorem}
  \label{theorem: unknot}
  The CAT unknotting problem for 2-spheres in  $\#_7 (S^2 \times S^2)$ is unsolvable.
\end{theorem}

Theorem \ref{theorem: unknot} generalizes to all surfaces, orientable and non-orientable; see Theorem \ref{theorem: sfceuns}.

A question that turns out to be related to Question \ref{question: 2knot} is the following.

\begin{question}
  \label{question: uswp}
Is there a 2-sphere in $S^4$ whose group has unsolvable word problem?
\end{question}

It is proved in [GGS] that a positive answer to Question \ref{question: uswp} implies a negative answer to Question \ref{question: 2knot}.

In [G3] we showed that there is an orientable surface of genus 18 in $S^4$ whose group has unsolvable word problem. Here we show 

\begin{theorem}
  \label{theorem: uswp}
  There is a 2-sphere in $\#_7 (S^2 \times S^2)$ whose group has unsolvable word problem.
\end{theorem}

In [G3] we showed that there is an orientable surface of genus 40 in $S^4$ whose group contains an isomorphic copy of every recursively presented group. Here we show

\begin{theorem}
  \label{theorem: univ}
  There is a 2-sphere in  $\#_{12} (S^2 \times S^2)$ whose group contains an isomorphic copy of every recursively presented group.
\end{theorem}

Again, Theorems \ref{theorem: uswp} and \ref{theorem: univ} generalize to all surfaces; see Theorems \ref{theorem: sfceuswp} and \ref{theorem: sfceuniv}.

$F(a,b,c,...)$ will denote the free group on $\{a,b,c,...\}$. We will write $\mathbb{Z}(a)$ for $F(a)$, and $\mathbb{Z}_2(x)$ for the cyclic group of order 2 generated by $x$.

In the context of group presentations, by a mild abuse of terminology we do not distinguish between relators and relations.

$[x,y]$ will denote $x^{-1}y^{-1}xy$, and $\langle\!\langle \; \rangle\!\rangle$ will denote normal closure.

If $F$ is a surface in a PL 4-manifold $X$, then the {\it group} of $F$, denoted by $\pi(F)$, will mean $\pi_1(X \setminus F)$.

$N(...)$ will denote regular neighborhood.

\section{PL Unrecognizability of $\#_7(S^2 \times S^2)$}

Let $P$ be a finite presentation $(x_1,...,x_m : r_1,...,r_n)$ of a group $G$. Let $V$ be the 5-manifold obtained by attaching $m$ orientable 1-handles to the 5-ball $B^5$. Thus $V$ is the boundary connected sum of $m$ copies of $S^1 \times D^4$, $\natural_{m}(S^1 \times D^4)$. Choosing a basepoint in $\partial B^5 \cap \partial V$, we have $\pi_1(\partial V) \cong \pi_1(V) \cong F(x_1,...,x_m)$, where the $x_i$ are represented by the cores of the 1-handles. If $\gamma$ is a simple closed curve in $\partial V$, then the choice of an arc $\alpha$ in $\partial V$ from the basepoint to $\gamma$, and an orientation of $\gamma$, determines an element $r \in F(x_1,...,x_m)$.  Let $\gamma_j$, $1 \le j \le n$, be disjoint simple closed curves in $\partial V$ such that, with some choice of arcs and orientations as above, $\gamma_j$ determines $r_j$, $1 \le j \le n$. There is an obvious embedding of $V$ in $\mathbb{R}^5$. Then $\gamma_j$ bounds a disk in $\mathbb{R}^5$, and we define the 0-framing of the normal bundle of $\gamma_j$ in $\partial V$ to be the one that extends over the normal bundle of the disk. Let $W(P)$ be the 5-manifold obtained by attaching 2-handles $h_j$ to $V$ along $\gamma_j$, $1 \le j \le n$, with the 0-framing. Clearly $\pi_1(W(P)) \cong G$. Let $X(P) = \partial W(P)$. Since inclusion $X(P) \to W(P)$ induces an isomorphism on $\pi_1$, we also have $\pi_1(X(P)) \cong G$.

{\it Andrews-Curtis (AC)-equivalence} is the equivalence relation on finite group presentations generated by the following transformations ({\it AC-moves}) on a presentation $(x_1,...,x_m : r_1,...,r_n)$ (see [AC]).

(1) replace some $r_j$ by $r_j^{-1}$

(2)  replace some $r_j$ by $u^{-1}r_ju$, where $u$ is a word in $x_1,...,x_m$

(3) replace some $r_j$ by $r_jr_k$, for some $k \ne j$

(4) add a generator $x_{m+1}$ and a relation $x_{m+1} = w$, where $w$ is a word in $x_1,...,x_m$.

It is easy to see that these moves do not alter the PL homeomorphism type of $W(P)$. More specifically, they correspond to:

(1) re-orienting $\gamma_j$

(2) choosing a different arc from the basepoint to $\gamma_j$

(3) sliding the 2-handle $h_j$ over $h_k$

(4) adding an extra 1-handle to $V$ and attaching a cancelling 2-handle.

Thus we have the following well-known theorem.
            
\begin{theorem}
 \label{theorem: AC} 
  If $P$ and $P'$ are AC-equivalent then $W(P) \overset{PL}{\cong} W(P')$. Hence $X(P) \overset{PL}{\cong} X(P')$.
\end{theorem}

Let $[1]_k$, $k \ge 0$, denote the $k$-term sequence $1,1,...,1$.

\begin{corollary}
  \label{corollary: ACtrivgroup}
  If a finite presentation $P$ is AC-equivalent to $(\emptyset : [1]_k)$ then $W(P) \overset{PL}{\cong} \natural_k(S^2 \times D^3)$. Hence $X(P) \overset{PL}{\cong} \#_k(S^2 \times S^2)$.
\end{corollary}

The following lemma is elementary.

\begin{lemma}
  \label{lemma: easy}
  Let $P = (x_1,...,x_m : x_1 = w, r_1,...,r_n)$, where $w$ is a word in $x_2,...,x_m$.  Then $P$ is AC-equivalent to $(x_2,...,x_m : r'_1,...,r'_n)$, where $r'_j$ is obtained from $r_j$ by replacing all occurrences of $x_1^{\pm1}$ by $w^{\pm1}$, $1 \le j \le n$.
  \end{lemma} 

Let $\Gamma$ be the group denoted by $\Gamma_6$ in [KLR]. Then $\Gamma$ has unsolvable word problem, and has a presentation $Q = (s,d,z: r_1,...,r_9)$. (Here $s$ is the element denoted by $s_1$ in [KLR].) 

It is shown in [KLR, Proposition 2.4] that $\langle\!\langle z \rangle\!\rangle = \Gamma$. In fact more is true.

  \begin{proposition}
\label{proposition: kill Gamma}
 $(s,d,z : r_1,...,r_9,z)$ is AC-equivalent to $(\emptyset : [1]_7)$.
\end{proposition}

\begin{proof}

In [KLR, Proposition 2.3], the group $\Gamma$ is initially defined by a presentation with 9 generators and 22 relators. However, it is shown that 7 of the relators are redundant; discarding these gives the following 9-generator 15-relator presentation $Q_1$:

   generators $s,d,z,t,c,y,e,k,s_2$, and relations

   (i)$'$    $d^4s = sd$

   (ii)$'$   $es = se^4$

   (iii)$'$  $sc = cs$

   (iv)$'$   $d^j F_j e^j c = c d^j E_j e^j$,  $j = 1,2,3$

   (v)$'$    $dt = td$

   (vi)$'$   $P^{-1} t P k = k P^{-1} t P$

   (vii)$'$  $dye = y$

   togther with

   $(\dagger)$  $t = z^{-1}dz, c = z^{-1}tz, y = z^{-1}cz, e = ydy^{-1}$, $k = y^{-1}ty$, and $s_2 = y^{-1} s^{-1} y$. 

    Here, $E_j$ and $F_j$, $j = 1,2,3,$ and $P$ are words in $s$ and $s_2$.

    Let $r'_1,...,r'_9$  be relators corresponding to the 9 relations (i)$'$ - (vii)$'$. Using the relations $(\dagger)$, $t,c,y,e$ and $k$ can be successively expressed as words in $d$ and $z$, and $s_2$ can be expressed as a word in $s,d$ and $z$. This shows that $\Gamma$ has presentation $Q = (s,d,z: r_1,...,r_9)$, where $r_1,...,r_9$ are the words in $s,d$ and $z$ obtained by making the above substitutions in $r'_1,...,r'_9$.

    In fact, $Q_1$ and $Q$ are AC-equivalent. To see this, we have $Q_1 =$

    $(s,d,z,t,c,y,e,k,s_2 :  t = z^{-1}dz, c = z^{-1}tz, y = z^{-1}cz, e = ydy^{-1}, k = y^{-1}ty, s_2 = y^{-1} s^{-1} y, r'_1,...,r'_9)$

    Let $T$ denote $z^{-1}dz$. Then, by Lemma \ref{lemma: easy}, $Q_1$ is AC-equivalent to

    $Q_2 = (s,d,z,c,y,e,k,s_2 :  c = z^{-1}Tz, y = z^{-1}cz, e = ydy^{-1}, k = y^{-1}Ty, s_2 = y^{-1} s^{-1} y, r^{(2)}_1,...,r^{(2)}_9)$

    where $r^{(2)}_j$ is obtained from $r'_j$ by replacing each occurrence of $t^{\pm 1}$ by $T^{\pm 1}$.

    Now let $C$ be the word $z^{-1}Tz$ in $d$ and $z$. Then $Q_2$ is AC-equivalent to

    $Q_3 = (s,d,z,y,e,k,s_2 :  y = z^{-1}Cz, e = ydy^{-1}, k = y^{-1}Ty, s_2 = y^{-1} s^{-1} y, r^{(3)}_1,...,r^{(3)}_9)$

where  $r^{(3)}_j$ is obtained from $r^{(2)}_j$ by replacing each occurrence of $c^{\pm 1}$ by $C^{\pm 1}$.
    
Continuing in this way, we eventually get that $Q_1$ is AC-equivalent to $(s,d,z,: r_1,...,r_9) = Q$.

It follows that $Q^+ = (s,d,z : r_1,...,r_9, z)$ is AC-equivalent to the presentation $Q_{1}^+$ obtained by adding the relator $z$ to $Q_1$. It is easy to see that $Q_{1}^+$ is AC-equivalent to 
$$(s,d,z,t,c,y,e,k,s_2 :  t = d, c = d, y = d, e = d, k = d, s_2 = d^{-1} s^{-1} d, r'_1,...,r'_9, z)$$

   Using $r'_3$, i.e. $sc = cs$, and the relation $c = d$, by an AC-equivalence we can replace the relation $s_2 = d^{-1} s^{-1} d$ by $s_2 = s^{-1}$.

  Using Lemma \ref{lemma: easy}, this new presentation is AC-equivalent to $R = (s,d,z : r''_1,...,r''_9, z)$, where $r''_j$ is obtained from $r'_j$ by making the substitutions  $t = d, c = d, y = d, e = d, k = d, s_2 = s^{-1}$.

   Since $r'_1,...,r'_9$ are given by (i)$'$ - (vii)$'$ above, we see that $r''_1,...,r''_9$ correspond to the relations

   (i)$''$    $d^4s = sd$

   (ii)$''$   $ds = sd^4$

   (iii)$''$  $sd = ds$

   (iv)$''$   $\bar{F}_j d = d \bar{E}_j$,  $j = 1,2,3$

   (v)$''$    $d^2 = d^2$

   (vi)$''$   $\bar{P}^{-1} d \bar{P} d = d \bar{P}^{-1} d \bar{P}$

   (vii)$''$  $d^3 = d$,

   where, if $S$ is a word in $s = s_1$ and $s_2$, then $\bar{S}$ denotes the result of substituting $s_2 =  s^{-1}$. Thus $\bar{S}$ is $s^{\epsilon_1 - \epsilon_2}$, where $\epsilon_i$ is the exponent sum of $s_i$ in $S$, $i = 1,2$.   

Note that (v)$''$ is equivalent to the relator $1$.  Also, by AC-equivalences, using (iii)$''$, (i)$''$ can be replaced by $d^3 = 1$, and then, using (vii)$''$, by $d = 1$. This can then be used to replace (iv)$''$ by $\bar{F}_j = \bar{E}_j$,  $j = 1,2,3$, and each of (ii)$''$, (iii)$''$, (vi)$''$, and (vii)$''$ by 1. So we get that $R$ is AC-equivalent to $R' = (s,d,z : d, [1]_5, \bar{F}_j = \bar{E}_j, j = 1,2,3, z)$. It is easy to see from [KLR, Appendix A] that $\bar{E}_1 = s^{-1}$ and $\bar{F}_1 = s^{-2}$, so $R'$ is AC-equivalent to $(s,d,z : s,d,z, [1]_7)$, and hence to $(\emptyset : [1]_7)$. 
 \end{proof}

\begin{lemma}
\label{lemma: w = 1}
Let $G$ be a group with presentation $(x_1,...,x_m : r_1,...,r_n)$. Let $P$ be a presentation of the form $(x_1,...,x_m,y_1,...,y_p : r_1,...,r_n, s_1,...,s_q)$, where each $s_i$ is a word in $x_1,...,x_m,y_1,..,y_p$. Suppose $s_1 = uwv$, where $w$ is a word in $x_1,...,x_m$ such that $w \overset{G}{=} 1$. Then $P$ is AC-equivalent to the presentation obtained from $P$ by replacing $s_1$ by $uv$.
\end{lemma}

\begin{proof}
Since $w \overset{G}{=} 1$, $w = \prod_{k=1}^l u_k^{-1}r^{\epsilon_k}_{j_k}u_k$ in $F(x_1,...,x_m)$, for some $u_k \in F(x_1,...,x_m)$ and $\epsilon_k = \pm1$, $1 \le k \le l$. We prove the lemma by induction on $l$, the case $l = 0$ being when $w$ is the empty word 1.

Perform the following sequence of AC-moves on the set of relators of $P$ (where $\hat{r}$ indicates that the relator $r$ has been removed). In each of these moves $s_2,...,s_q$ are unchanged, so to simplify the notation we omit them from the description of the moves.

If $\epsilon_1 = -1$:    $\{r_1,...,r_n,s_1\} \to \{r_1,...,r_n,wvu\} \to \{r_1,...,\hat{r}_{j_1},...,r_n, u_1^{-1}r_{j_1}u_1, wvu\}$

      $\to \{r_1,...,\hat{r}_{j_1},...,r_n, u_1^{-1}r_{j_1}u_1, u_1^{-1}r_{j_1}u_1wvu\} = \{r_1,...,\hat{r}_{j_1},...,r_n, u_1^{-1}r_{j_1}u_1, s'_1\} \to \{r_1,...,r_n, s'_1\}$

If $\epsilon_1 = +1$:   before the first move above replace $r_{j_1}$ by $r^{-1}_{j_1}$, and after the last move replace $r^{-1}_{j_1}$ by $r_{j_1}$.

Note that in both cases $s'_1 = w'vu$, where $w' = \prod_{k=2}^l u_k^{-1}r^{\epsilon_k}_{j_k}u_k$. Hence the result follows by induction on $l$.
\end{proof}

Recall the presentation $Q = (s,d,z: r_1,...,r_9)$ of $\Gamma$. Let $w$ be a word in $s,d,z$, and let $G_w$ be the group with the following presentation $P_w$:

generators $s,d,z,a,b,\beta$, and relations $r_1,...,r_9$ together with

(i) $a^{-1}za = \beta^{-1} b \beta$

(ii) $a^{-2}za^2 = b^{-1} \beta b$

(iii) $z^{-1}az = b^{-2} \beta b^2$

(iv) $a [a^{-1},w] a^{-1} = b^{-3} \beta b^3$

\begin{proposition}
  \label{proposition: Pw}
  (1) If $w \overset{\Gamma}{=} 1$ then $P_w$ is AC-equivalent to $(\emptyset : [1]_7)$.

  (2) If $w \overset{\Gamma}{\neq} 1$ then $G_w \neq 1$.
\end{proposition} 

\begin{proof}

(1) Suppose $w \overset{\Gamma}{=} 1$. Then by Lemma \ref{lemma: w = 1}, by an AC-equivalence we can replace relation (iv) by $\beta = 1$. Then, by further AC-equivalences, (iii) can be replaced by $a = 1$, (ii) by $z = 1$, and then (i) by $b = 1$. Hence $P_w$ is AC-equivalent to $(s,d,z,a,b,\beta : a,b,\beta, r_1,...,r_9, z)$, which is AC-equivalent to  $Q^+ = (s,d,z : r_1,...,r_9, z)$. The result now follows from Proposition \ref{proposition: kill Gamma}. 

(2) We claim that the elements on the right-hand side of relations (i) - (iv) are a basis for a free subgroup of $F(b,\beta)$. To see this, let $Y$ denote the set of words listed on the right-hand side of the relations. By examining the possible cancellation in $y_1^{\epsilon_{1}}y_2^{\epsilon_{2}}$, where $y_1$ and $y_2$ are distinct elements of $Y$ and $\epsilon_i = \pm1, i = 1,2$, it is easy to see that a non-empty reduced word in the elements of $Y$ is non-trivial when expressed as a reduced word in $b$ and $\beta$. Thus $Y$ is a basis for a free subgroup of the free group $F(b,\beta)$. 

Next note that since the exponent sum of $z$ in each relator of $Q_1$ is zero, $z$ has infinite order in $\Gamma$. It follows (arguing as in the previous paragraph) that if $w \overset{\Gamma}{\neq} 1$ then the elements on the left-hand side of (i) - (iv) are a basis for a free subgroup of $\Gamma * \mathbb Z(a)$. Hence if $w \overset{\Gamma}{\neq} 1$ then $G_w$ is a free product with amalgamation $(\Gamma * \mathbb Z(a))*_{E}F(b,\beta)$, where $E$ is a free group of rank 4.
\end{proof}

\begin{proof}[Proof of Theorem \ref{theorem: rec7}]
 If $w \overset{\Gamma}{=} 1$ then $P_w$ is AC-equivalent to $(\emptyset : [1]_7)$ by Proposition \ref{proposition: Pw}(1), and hence $X(P_w) \overset{PL}{\cong} \#_{7}(S^2 \times S^2)$ by Theorem \ref{theorem: AC}. On the other hand, if $w \overset{\Gamma}{\neq} 1$ then $\pi_1(X(P_w)) \cong G_w \neq 1$ by Proposition \ref{proposition: Pw}(2). Therefore $X(P_w) \overset{PL}{\cong} \#_{7}(S^2 \times S^2)$ if and only if $w \overset{\Gamma}{=} 1$. Since $\Gamma$ has unsolvable word problem, the theorem is proved. 
\end{proof}

\begin{remark}
The basic structure of all the unrecognizability proofs we have mentioned is the following. Start with a finite presentation of a group $G$ with unsolvable word problem and a word $w$ in the generators. Next construct a presentation $P_w$ of a group $G_w$ such that if $w \overset{G}{\neq} 1$ then $G_w \neq 1$, and if $w \overset{G}{=} 1$ then the boundary $X(P_w)$ of the 5-manifold $W(P_w)$ is homeomorphic to $\#_k(S^2 \times S^2)$, either in the PL or TOP category. This requires choosing the 0-framing for the attaching maps of the 2-handles in the construction of $W(P_w)$, and this can be done algorithmically, as explained in detail in [S1]; see also [FHK]. However, we point out the following alternative aproach. For the attaching circle of each 2-handle first choose any one of the two possible framings; this then determines the other framing. Hence one can algorithmically construct $2^n$ 5-manifolds $W_i(P_w)$, with boundaries $X_i(P_w)$, $1 \le i \le 2^n$. Then $w \overset{G}{=} 1$ if and only if some $X_i(P_w)$ is (PL or TOP) homeomorphic to $\#_k(S^2 \times S^2)$. But if the corresponding recognition problem for $\#_k(S^2 \times S^2)$ were solvable then the latter property would be decidable.
\end{remark}

\section{Unknotting surfaces in 4-manifolds}

We now turn to the unknotting problem for surfaces in 4-manifolds. By a {\it surface} we will mean a closed, connected 2-manifold, and if $X$ is a compact, connected PL 4-manifold, by a {\it surface in $X$} we will mean a locally flat PL embedded surface in Int $X$.

To define what it means for a surface in a 4-manifold to be unknotted, first recall that in $S^4$, an orientable surface is CAT {\it unknotted} if it bounds a CAT handlebody, and a non-orientable surface is CAT {\it unknotted} if there is a CAT isotopy of $S^4$ taking it to a connected sum of copies of the standard projective planes $P_{\pm}$; see [HK]. Next, a surface $F$ in a PL 4-ball $B$ is CAT {\it unknotted} if for some (and hence any) PL embedding $e: B \to S^4$, $e(F)$ is CAT unknotted in $S^4$. Finally, a surface $F$ in a compact, connected PL 4-manifold $X$ is CAT {\it unknotted} if for some (and hence any) PL 4-ball $B \subset$ Int $X$, there is a CAT isotopy of $X$ taking $F$ to a CAT unknotted surface in $B$. In particular, a 2-sphere in $X$ is CAT unknotted if and only if it bounds a CAT 3-ball in $X$.

Note that if a surface $F$ in $X$ is CAT unknotted then $\pi(F) \cong C*\pi_1(X)$, where $C \cong \mathbb {Z}$ or $\mathbb {Z}_2$, according as $F$ is orientable or non-orientable.

We say that the CAT {\it unknotting problem for surfaces of a given homeomorphism type in $X$ is solvable} if there is an algorithm to decide whether or not a given surface in $X$ of that homeomorphism type is CAT unknotted.

The key cases are when the surface is either $S^2$ or $\mathbb{R}P^2$, and for these we have the following proposition. For part (1) in the the case $d = -1$ (where $\#_0(S^2 \times S^2)$ is to be interpreted as $S^4$), see [L2].

\begin{proposition}
  \label{proposition: killerg}
  Let $P = (x_1,...,x_m : r_1,...,r_n)$ be a presentation of a group $G$, with deficiency $d = n - m \ge -1$. Suppose there exists $g \in G$ such that the presentation $P^+ = (x_1,...,x_m : r_1,...,r_n, g)$ is AC-equivalent to $(\emptyset : [1]_{d+1})$. 

  (1) There is a 2-sphere $S$ in $\#_{d+1}(S^2 \times S^2)$ such that $\pi(S) \cong G$, where $g$ is represented by a meridian of $S$.

  (2) If $P$ is AC-equivalent to $(g : [1]_{d+1})$ then $S$ is PL unknotted.

  (3) There is a projective plane $K$ in $\#_{d+1}(S^2 \times S^2)$ such that $G$ embeds in $\pi(K)$. 
    \end{proposition}

\begin{proof}[Proof of Proposition \ref{proposition: killerg}(1)] 

  Let $B^3 \subset B^5$ be the standard inclusion, with boundaries $S^2 \subset S^4$. Let $V$ be the result of attaching $m$ orientable 1-handles to $B^5$, with attaching maps disjoint from $N(S^2)$. Then $B^3$ is a properly embedded 3-ball in $V$, with boundary $S^2$ in $\partial V$. We identify $\pi_1(V)$ with the free group $F(x_1,...,x_m)$. Then $\pi_1(V \setminus$ int$N(B^3)) \cong F(x_1,...,x_m,y)$, where $y$ is represented by a meridian of $B^3$.

  Consider the 5-manifold $Y$ obtained by attaching $(n+1)$ 2-handles to $V \setminus$ int$N(B^3)$, with attaching maps in $\partial V \setminus$ int$N(S^2)$ representing $r_1,...,r_n$ and $gy^{-1}$ in $F(x_1,...,x_m,y)$. Thus $Y \overset{PL} \cong W(P')$, where $P'$ is the presentation $(x_1,...,x_m,y : r_1,...,r_n, gy^{-1})$. Clearly $P'$ is AC-equivalent to $P$, hence $Y \overset{PL} \cong W(P)$  and $\pi_1(Y) \cong G$. Let $W = Y \cup N(B^3)$. Since $y$ is a meridian of $B^3$, $W \overset{PL} \cong W(P'^+)$, where $P'^+ = (x_1,...,x_m,y : r_1,...,r_n,gy^{-1},y)$. Clearly $P'^+$ is AC-equivalent to $P^+$. Therefore $W(P'^+) \overset{PL} \cong W(P^+) \overset{PL} \cong \natural_{d+1}(S^2 \times D^3)$ by Corollary \ref{corollary: ACtrivgroup} and our hypothesis on $P^+$. Then $B^3$ in $W$ corresponds to a 3-ball $B$ in $\natural_{d+1}(S^2 \times D^3)$ with $\pi(B) \cong \pi_1(Y) \cong G$, where $g$ is a meridian of $B$. The boundary of $B$ is a 2-sphere $S$ in $\#_{d+1}(S^2 \times S^2)$ such that $\pi(S) \cong G$, where $g$ is a meridian of $S$.
\end{proof}

Before proving part (2), we make the following definition. Let $M$ be a PL 4-manifold and $B \subset M$ a properly embedded locally flat PL 3-ball. We say that $B$ is {\it PL unknotted} in $M$ if there exists a locally flat PL 4-ball $\Delta \subset M$ such that $\partial \Delta = B \cup_{\partial} B_0$, where $B_0 = \Delta \cap \partial M$ is a 3-ball. Note that $S = \partial B$ bounds $B_0 \subset \partial M$ and hence is PL unknotted in $\partial M$.

\begin{proof}[Proof of Proposition \ref{proposition: killerg}(2)] 

  Recall from the proof of part (1) that we have $(\natural_{d+1}(S^2 \times D^3), B) \overset{PL} \cong (W,B^3)$, where $W = Y \cup N(B^3)$. Suppose that $P$ is AC-equivalent to $(g : [1]_{d+1})$; we will show that $B^3$ is PL unknotted in $W$. This implies that $B$ is PL unknotted in $(\natural_{d+1}(S^2 \times D^3)$, and hence that $S$ is PL unknotted in $\#_{d+1}(S^2 \times S^2)$. 

  Since $Y \overset{PL} \cong W(P)$ from the proof of part (1) above, the hypothesis on $P$ implies that there is a PL homeomorphism $f : \natural_{d+1}(S^2 \times D^3) \natural(S^1 \times D^4) \to Y$. We have $\partial D^4 = S^3 = B_{+}^3 \cup B_{-}^3$, the union of its upper and lower hemispheres. Then we may assume that the 4-ball in $S^1 \times S^3$ realizing the boundary connected sum of $S^1 \times D^4$ with $\natural_{d+1}(S^2 \times D^3)$ lies in $J \times$ int$B_{-}^3$ for some interval $J \subset S^1$. 

  We have $N(B^3) = D^2 \times B^3$, and $\partial N(B^3) = (S^1 \times B^3) \cup_{\partial} (D^2 \times S^2)$, where $S^1 \times B^3 = N(B^3) \cap \partial Y$. Since homotopy implies PL isotopy for PL simple closed curves in PL 4-manifolds, and by uniqueness of regular neighborhoods, we may assume that  $f(S^1 \times B_{+}^3) = S^1 \times B^3$. We may also assume that for some $p \in S^1 \setminus J$, $f(\{p\} \times B_{+}^3) = \{p'\} \times B^3$ for some $p' \in S^1$. Let $D$ be the 4-ball $f(\{p\} \times D^4)$, properly embedded in $Y$. Then $\partial D = (\{p'\} \times B^3) \cup_{\partial} B_0$, where $B_0 = D \cap \partial W$ is the 3-ball $f(\{p\} \times B_{-}^3)$.

  Let $I$ be the radial interval in $D^2$ with $\partial I = \{0, p'\}$, and consider $I \times B^3 \subset D^2 \times B^3 = N(B^3)$. Then $\Delta = D \cup_{\{p'\} \times B^3} (I \times B^3)$ is a 4-ball in $W$ that demonstrates the PL unknottedness of $B^3 = B^3 \times \{0\}$ in $W$.
\end{proof}

We wish to prove Theorem \ref{theorem: unknot} and its generalization to all surfaces, Theorem \ref{theorem: sfceuns}. To deal with non-orientable surfaces we will need part (3) of Proposition \ref{proposition: killerg}, so the next construction is done with that in mind.

Recalling the notation in Proposition \ref{proposition: killerg}, let $R$ be the presentation with generators $x_1,...,x_m, a, \alpha, b, \beta$, and relators $r_1,...,r_n$ together with

(i) $\alpha^{-1} a \alpha = \beta b \beta^{-1}$

(ii) $a^{-1} \alpha a = b^{-1} \beta b$

(iii) $\alpha^{-2} a^{-1} (g \alpha) a \alpha^2 = \beta^{-1} (b^{-1} \beta b) \beta$

(iv) $\alpha^{-3} [h,a] \alpha^3 = \beta^{-2} (b^{-1} \beta b) \beta^2$

where $h$ is any element in $G \setminus\{1\}$.

Let $R^+$ be the presentation obtained from $R$ by adding the relation $a = 1$.

\begin{lemma}
  \label{lemma: R+}
  $R^+$ is AC-equivalent to $(\emptyset : [1]_{d+1})$.
\end{lemma}

\begin{proof}
Using the relation $a = 1$, relation (i) can be replaced by $b = 1$, and (iv) can be replaced by $\beta = 1$. From this (ii) can be replaced by $\alpha = 1$ and (iii) by $g = 1$. It is easy to see that these can be realized by AC-equivalences, and hence $R^+$ is AC-equivalent to $P^+$, which is AC-equivalent to $(\emptyset: [1]_{d+1})$ by hypothesis.
\end{proof}

Let $H$ be the group presented by $R$. Then Lemma \ref{lemma: R+} and Proposition \ref{proposition: killerg}(1) imply

\begin{corollary}
  \label{corollary: 2sphere} 
 There is a 2-sphere $S$ in $\#_{d+1}(S^2 \times S^2)$ such that $\pi(S) \cong H$, where $a$ is represented by a meridian of $S$.
\end{corollary}

\begin{proof}[Proof of Proposition \ref{proposition: killerg}(3)]
Let $S$ be the 2-sphere in Corollary \ref{corollary: 2sphere}, and let $K = S\#P_{\pm}$, where $P_{\pm}$ is any unknotted projective plane in $S^4$. Then $\pi(K) \cong H/\langle\!\langle a^2 \rangle\!\rangle$.

 Consider the presentation obtained by adding the relation $a^2 = 1$ to $R$. By (i), this implies $b^2 = 1$, so we can now regard the right-hand side of relations (i) - (iv) as elements of $\mathbb{Z}_2(b) * \mathbb{Z}(\beta)$ and the left-hand side as elements of $G * \mathbb{Z}_2(a) * \mathbb{Z}(\alpha)$. It is easy to see that each of these sets of elements generates a subgroup isomorphic to $\mathbb{Z}_2 * F_3$.
 Formally, let $\mathbb{Z}_2(x)$ be the group $\langle x: x^2 = 1 \rangle$, and let $E$ the free product $\mathbb{Z}_2(x) * F(y_1,y_2,y_3)$. Let $\rho: E \to \mathbb{Z}_2(b)*\mathbb{Z}(\beta)$ be the homomorphism defined by

  $x \mapsto  \beta b \beta^{-1}$

  $y_1 \mapsto  b^{-1} \beta b$

  $y_2 \mapsto  \beta^{-1} (b^{-1} \beta b) \beta$

  $y_3 \mapsto  \beta^{-2} (b^{-1} \beta b) \beta^2$

Similarly, let $\lambda : E \to G *\mathbb{Z}_2(a) * \mathbb{Z}(\alpha)$ be the homomorphism defined by

  $x \mapsto  \alpha^{-1} a \alpha$
 
  $y_1 \mapsto  a^{-1} \alpha a$

  $y_2 \mapsto  \alpha^{-2} a^{-1} (g \alpha) a \alpha^2$ 

  $y_3 \mapsto  \alpha^{-3} [h,a] \alpha^3$

It is straightforward to check that if $v$ is a non-trivial word in $\{x,y_1,y_2,y_3\}$ that is reduced as a word in the free product $\mathbb{Z}_2(x) * \mathbb{Z}(y_1)*\mathbb{Z}(y_2)*\mathbb{Z}(y_3)$, then $\rho(v)$ and $\lambda(v)$ can be expressed as non-trivial reduced words in the free products $\mathbb{Z}_2(b) *\mathbb{Z}(\beta)$ and $G *\mathbb{Z}_2(a) *\mathbb{Z}(\alpha)$, respectively. Hence $\rho$ and $\lambda$ are injective.

It follows that $H/\langle\!\langle a^2 \rangle\!\rangle$ is a free product with amalgamation $(G * \mathbb{Z}_2 * \mathbb{Z}) *_{\mathbb{Z}_2 * F_3}(\mathbb{Z}_2 * \mathbb{Z})$. Therefore $G$ embeds in $H/\langle\!\langle a^2 \rangle\!\rangle \cong \pi(K)$.
\end{proof}

\begin{remark}
  \label{remark: inforder}
  In the presentation $R$, if $g$ has infinite order in $G$ then the left-hand side of relation (iii) can be replaced by  $\alpha^{-2} g \alpha^2$. Note that since $\langle\!\langle g \rangle\!\rangle = G$, this is the case if $H_1(G) \cong \mathbb{Z}$.
\end{remark}
  
We now proceed towards the proof of Theorem \ref{theorem: unknot} and its generalization to all surfaces, Theorem \ref{theorem: sfceuns}. We start with the presentation $Q = (s,d,z : r_1,...,r_9)$ of the group $\Gamma$ with unsolvable word problem [KLR]. 

Let $w$ be a word in $s,d,z$. Let $Q_w$ be the presentation with generators $s,d,z,a,\alpha,b,\beta$, and relators $r_1,...,r_9$ together with relations (i) - (iv) in the definition of $R$ above, in which $g = z$ and $h = w$. (Here, we allow the possibility that $w \overset{\Gamma}{=} 1$. Note also that since $z$ has infinite order in $\Gamma$ we can take the left-hand side of relation (iii) to be $\alpha^{-2} z \alpha^2$; see Remark \ref{remark: inforder}.)

Let $H_w$ be the group presented by $Q_w$.

\begin{proposition}
  \label{proposition: Hw}
  
  (1) If $w \overset{\Gamma}{=} 1$ then $Q_w$ is AC-equivalent to $(a : [1]_7)$.

  (2) If $w \overset{\Gamma}{\neq} 1$ then $H_w \ncong \mathbb{Z}$.

  (3) If $w \overset{\Gamma}{\neq} 1$ then $H_w /\langle\!\langle a^2 \rangle\!\rangle \ncong \mathbb{Z}_2$.
\end{proposition}

\begin{proof}

  (1) Suppose $w \overset{\Gamma}{=} 1$. Then by Lemma \ref{lemma: w = 1}, by an AC-equivalence relation (iv) can be replaced by $\beta = 1$. Then (ii) can be replaced by $\alpha = 1$ and  (iii) by $z = 1$. Finally, (i) can be replaced by $a = b$. Thus $Q_w$ is AC-equivalent to $(s,d,z,a : r_1,...,r_9, z)$. The result now follows from Proposition \ref{proposition: kill Gamma}. 

  (3) The proof of Proposition \ref{proposition: killerg}(3) shows that $H_w/\langle\!\langle a^2 \rangle\!\rangle$ is a free product with amalgamation $(\Gamma * \mathbb{Z}_2 * \mathbb{Z}) *_{\mathbb{Z}_2 * F_3}(\mathbb{Z}_2 * \mathbb{Z})$.

  (2) Since $H_1(H_w) \cong \mathbb{Z}(a)$, this follows from (3).  
\end{proof}

Let $Q^{+}_w$ be the presentation obtained from $Q_w$ by adding the relation $a = 1$. 

\begin{lemma}
\label{lemma: Qw+}
$Q^{+}_w$ is AC-equivalent to $(\emptyset : [1]_7)$.
\end{lemma}

\begin{proof}
In the presentation $Q^{+}_w$, using the relation $a = 1$ relation (i) can be replaced by $b = 1$, and (iv) by $\beta = 1$. Then (ii) can be replaced by $\alpha = 1$, and (iii) by $z = 1$. Hence $Q^{+}_w$ is AC-equivalent to $(s,d,z :r_1,...,r_9,z)$. The result now follows from Proposition \ref{proposition: kill Gamma}.
\end{proof}

By Lemma \ref{lemma: Qw+}, taking $P = Q_w$, $G = H_w$, and $g = a$ in Proposition \ref{proposition: killerg}, part (1) of that proposition gives a 2-sphere $S_w$ in $\#_{7}(S^2 \times S^2)$ with $\pi(S_w) \cong H_w$ and meridian $a$.

\begin{proof}[Proof of Theorem \ref{theorem: unknot}]

   If $w \overset{\Gamma}{\neq} 1$ then $\pi(S_w) \cong H_w \ncong \mathbb{Z}$ by Proposition \ref{proposition: Hw}(2). Hence $S_w$ is CAT knotted.

   If $w \overset{\Gamma}{=} 1$ then $Q_w$ is AC-equivalent to $(a : [1]_7)$ by Proposition \ref{proposition: Hw}(1). Therefore $S_w$ is CAT unknotted by Proposition \ref{proposition: killerg}(2).

   Since $\Gamma$ has unsolvable word problem, this completes the proof.
\end{proof}

\begin{theorem}
\label{theorem: sfceuns}
The CAT unknotting problem for surfaces of any homeomorphism type in $\#_{7}(S^2 \times S^2)$ is unsolvable.
\end{theorem}

\begin{proof}
For a given homeomorphism type, let $F \subset S^4$ be a PL unknotted surface of that type, and let $F_w = S_w \# F \subset \#_{7}(S^2 \times S^2)$. Note that $\pi(F) \cong \mathbb{Z}$ or $\mathbb{Z}_2$, according as $F$ is orientable or non-orientable.

If $w \overset{\Gamma}{=} 1$ then $S_w$ is CAT unknotted, and hence $F_w$ is CAT unknotted.

Suppose $w \overset{\Gamma}{\neq} 1$. If $F$ is orientable then $\pi(F_w) \cong \pi(S_w) \ncong \mathbb{Z}$, so $F_w$ is CAT knotted. If $F$ is non-orientable then $\pi(F_w) \cong \pi(S_w)/\langle\!\langle a^2 \rangle\!\rangle$. By Proposition \ref{proposition: Hw}(3) this is not $\mathbb{Z}_2$, so $F_w$ is CAT knotted.
\end{proof}

\begin{proof}[Proof of Theorem \ref{theorem: uswp}]

By Propositions \ref{proposition: kill Gamma} and \ref{proposition: killerg}(1) there is a 2-sphere $S$ in $\#_{7}(S^2 \times S^2)$ such that $\pi(S) \cong \Gamma$. 
\end{proof}

\begin{theorem}
\label{theorem: sfceuswp}
There is a surface of any given homeomorphism type in $\#_{7}(S^2 \times S^2)$ whose group has unsolvable word problem.
\end{theorem}

\begin{proof}
Taking $(G,g) = (\Gamma,z)$ in Proposition \ref{proposition: killerg} and using Proposition \ref{proposition: kill Gamma} we get a 2-sphere $S$ in $\#_{7}(S^2 \times S^2)$ and a projective plane $K$ in $\#_{7}(S^2 \times S^2)$ such that $\pi(S) \cong \Gamma$ and $\pi(K)$ has a subgroup isomorphic to $\Gamma$. This proves the theorem for 2-spheres and projective planes. 

For other surfaces, take the connected sum of $S$ (resp. $K$) with a PL unknotted orientable (resp. non-orientable) surface in $S^4$. The new surface has group $\pi(S)$ (resp. $\pi(K)$).
\end{proof}

A group $G$ is {\it universal} if it contains as a subgroup an isomorphic copy of every recursively presented group. The existence of a finitely presented universal group was proved by Higman [Hi]. Later, Valiev constructed a universal group having a presentation with 5 generators and 21 relators [V]. Recently, this was improved by Ren, who constructed a universal group with a 2-generator 13-relator presentation [Re].

\begin{theorem}
 \label{theorem: sfceuniv} 
There is a surface of any given homeomorphism type in $ \#_{12}(S^2 \times S^2)$ whose group is universal.
\end{theorem}

Let $U$ be the universal group constructed in [Re]. Then $U$ has a presentation $T = (c,r :  r_1,...,r_{13})$.

\begin{proposition}
    \label{proposition: Utriv}
    The presentation $T^+ = (c,r : r_1,...,r_{13},r)$ is AC-equivalent to $(\emptyset : [1]_{12})$.
  \end{proposition}

  \begin{proof}
    By Lemma \ref{lemma: easy}, $T^+$ is AC-equivalent to $(c : r'_1,...,r'_{13})$, where $r'_j$ is obtained from $r_j$ by deleting all occurrences of $r^{\pm 1}$. Thus each $r'_j$ is a power of $c$. On the other hand, we see from [Re, Section 7.1] that (for example) $r'_7 = c$. Hence $T^+$ is AC-equivalent to $(c : c, [1]_{12})$, which is AC-equivalent to $(\emptyset : [1]_{12})$.
\end{proof}

\begin{proof}[Proof of Theorem \ref{theorem: sfceuniv}]
By Proposition \ref{proposition: Utriv}, $T^+$ satisfies the condition in the hypothesis of Proposition \ref{proposition: killerg}. Hence, arguing as in the proof of Theorem \ref{theorem: sfceuswp}, for every homeomorphism type of surface, we get a surface $F$ of that type in $ \#_{12}(S^2 \times S^2)$ such that $U$ embeds in $\pi(F)$.
\end{proof}

Theorem \ref{theorem: univ} is a special case of Theorem \ref{theorem: sfceuniv}.

 Finally we point out that Theorems \ref{theorem: rec7} and \ref{theorem: sfceuns} hold for a more general class of 4-manifolds. For the recognition problem this extension is well known; see [CL, Section 5] for example.

 \begin{theorem}
 \label{theorem: arbmfld}  
    Let $X_0$ be a compact, connected PL 4-manifold.

    (1) The CAT recognition problem for $X_0 \# (\#_7(S^2 \times S^2))$ is unsolvable.

    (2) The CAT unknotting problem for surfaces of any homeomorphism type in $X_0 \# (\#_7(S^2 \times S^2))$ is unsolvable.
  \end{theorem}

  \begin{proof}

    (1) Let $P_w$ be as in Proposition \ref{proposition: Pw}. Define $X_w = X_0 \# X(P_w)$.

    If $w \overset{\Gamma}{=} 1$ then $X_w \overset{PL}{\cong} X_0 \# (\#_7(S^2 \times S^2))$.

    If $w \overset{\Gamma}{\neq} 1$ then $\pi_1(X_w) \cong \pi_1(X_0) * G_w$, where $G_w \neq 1$. Therefore, by uniqueness of free product decompositions of finitely generated groups,
    $\pi_1(X_w) \ncong \pi_1(X_0) \cong \pi_1(X_0 \# (\#_7(S^2 \times S^2)))$. Hence $X_w \overset{TOP}{\ncong} X_0 \# (\#_7(S^2 \times S^2))$.

 (2) Let $F_w$ be the surface in $\#_7 (S^2 \times S^2)$ in the proof of Theorem \ref{theorem: sfceuns} and $F'_w$ the corresponding surface in $X_0 \# (\#_7(S^2 \times S^2))$. 

 If $w \overset{\Gamma}{=} 1$ then $F_w$ is PL unknotted in $\#_7 (S^2 \times S^2)$, and hence $F'_w$ is PL unknotted in $X_0 \# (\#_7(S^2 \times S^2))$.

 If $w \overset{\Gamma}{\neq} 1$ then $\pi(F_w)$ is non-cyclic. Also, $\pi(F'_w) \cong \pi(F_w) * \pi_1(X_0)$. But the group of a TOP unknotted surface in $X_0 \# (\#_7(S^2 \times S^2))$ is isomorphic to $C*\pi_1(X_0)$, where $C$ is  $\mathbb Z$ or $\mathbb {Z}_2$, according as the surface is orientable or non-orientable. Hence, by the uniqueness of free product decompositions of finitely generated groups, $F'_w$ is TOP knotted in $X_0 \# (\#_7(S^2 \times S^2))$.
\end{proof}

\begin{remark}
 Part (2) of Theorem \ref{theorem: arbmfld} is only interesting for surfaces of relatively low genus. For it is shown in [G3] that in any compact, connected PL 4-manifold, the TOP (resp. PL) unknotting problem is unsolvable for orientable surfaces of genus $\ge 18$ (resp. $19$), and for non-orientable surfaces of genus $\ge 55$ (resp. $59$). (Here, the {\it genus} of a connected sum of $h$ projective planes is $h$.)  
\end{remark}

\medskip

{\bf Tool Use.} No sophisticated tools, including AI, were used in the preparation of this work.

\end{document}